\documentclass{amsart}

\usepackage{amsmath,amssymb,amsthm,latexsym}

\newtheorem{theorem}{Theorem}
\newcommand{\bt}{\begin{theorem}}
\newcommand{\et}{\end{theorem}}
\newtheorem{lemma}{Lemma}
\newcommand{\bl}{\begin{lemma}}
\newcommand{\el}{\end{lemma}}
\newtheorem{corollary}{Corollary}
\newcommand{\bc}{\begin{corollary}}
\newcommand{\ec}{\end{corollary}}
\newcommand{\bconj}{\begin{conjecture}}
\newcommand{\econj}{\end{conjecture}}
\newtheorem{problem}{Problem}
\newcommand{\bprob}{\begin{problem}}
\newcommand{\eprob}{\end{problem}}
\newcommand{\beq}{\begin{equation}}
\newcommand{\eeq}{\end{equation}}
\newcommand{\benum}{\begin{enumerate}}
\newcommand{\eenum}{\end{enumerate}}
\newcommand{\N}{\ensuremath{ \mathbf N }}
\newcommand{\Z}{\ensuremath{\mathbf Z}}

\newcommand{\mch}{\ensuremath{ \mathcal H}}

\newcommand{\bmat}{\left(\begin{matrix}}
\newcommand{\emat}{\end{matrix}\right)}
\newcommand{\bsmallmat}{\left(\begin{smallmatrix}}
\newcommand{\esmallmat}{\end{smallmatrix}\right)}

\title{Arithmetical structure of sumset intersections}

\author{Diego Marques}
\address{Departamento de Matem\' atica, Universidade de Bras\' ilia, Bras\' ilia, DF, Brazil}
\email{diego@mat.unb.br}
\author{Melvyn B.  Nathanson}
\address{Department of Mathematics\\Lehman College (CUNY)\\Bronx, NY 10468}
\email{melvyn.nathanson@lehman.cuny.edu}

\date{\today}


\thanks{D.M. would like to acknowledge the financial support provided by the National
Council for Scientific and Technological Development (CNPq).} 

\thanks{M.B.N.  supported in part by  PSC-CUNY Research Award Program grant 66197-00 54.}

\begin{document}

\begin{abstract}
The $h$-fold sumset of a set $A$ of integers is the set of all sums of $h$ not necessarily 
distinct elements of $A$.  
Let  $(A_q)_{q=1}^{\infty}$ be a strictly decreasing sequence of sets of integers and let 
$A = \bigcap_{q=1}^{\infty} A_q$.  Then $hA \subseteq   \bigcap_{q=1}^{\infty} hA_q$ 
for all $h \geq 1$.  Let $\mch(A_q) = \{h \geq 1: hA =  \bigcap_{q=1}^{\infty} hA_q\}$. 
The arithmetical structure of the sets $\mch(A_q)$ is unknown.  
It is proved that for every $h_0 \geq 2$ 
there exist sequences  $(A_q)_{q=1}^{\infty}$ such that 
$\{1,\ldots, h_0-1\} \subseteq \mch(A_q)$ but $h_0 \notin \mch(A_q)$ 
and also that there exist sequences  $(A_q)_{q=1}^{\infty}$ such that 
$\{1,  h_0 \} \subseteq \mch(A_q)$ but $\{2,3, \ldots, h_0-1\} \cap \mch(A_q) = \emptyset$.
\end{abstract}

\date{\today}

\maketitle 
\section{A sumset intersection problem} 
Let $\N = \{1,2,3,\ldots\}$ be the set of positive integers. 
A decreasing sequence of sets $(A_q)_{q=1}^{\infty}$ is 
\emph{strictly decreasing} if $A_q \neq A_{q+1}$ for all $q \in \N$ 
and \emph{asymptotically strictly decreasing} if $A_q \neq A_{q+1}$ 
for infinitely many $q \in \N$. 
A decreasing sequence of sets $(A_q)_{q=1}^{\infty}$ that is not asymptotically 
strictly decreasing is eventually constant. 

Let $X$ be an additive abelian semigroup and let $A$ be a subset of $X$.  
The $h$-fold sumset of $A$, denoted $hA$, is the set of all sums of $h$ not necessarily 
distinct elements of $A$.   
Nathanson~\cite{nath26Aus} posed the following problem:

\begin{quotation}
Let $(A_q)_{q=1}^{\infty}$ be a strictly decreasing sequence of sets 
in an additive abelian semigroup and let 
$A = \bigcap_{q=1}^{\infty}A_q$.  
Describe the set 
\[
\mch(A_q) = \left\{h \in \N: hA = \bigcap_{q=1}^{\infty} hA_q \right\}.
\] 
What sets of integers arise in this way? 
\end{quotation} 

Of course, $A = \bigcap_{q=1}^{\infty}A_q$ implies 
$1 \in \mch(A_q)$ for all sequences $(A_q)_{q=1}^{\infty}$. 

For all $q \in \N$, we have $A \subseteq A_q$ and so $hA \subseteq hA_q$.
Therefore, 
\[
hA \subseteq \bigcap_{q=1}^{\infty} hA_q. 
\]
Thus, $h \in \mch(A_q)$ if and only if $x \in hA_q$ for all $q \in \N$ implies $x \in hA$. 

If the decreasing sequence $(A_q)_{q=1}^{\infty}$ is eventually constant,
 then there exists $q_0 \in \N$ such that  
\[
A = \bigcap_{q=1}^{\infty} A_q = A_{q_0}
\]
and 
\[
hA = hA_{q_0} = \bigcap_{q=1}^{\infty} hA_q 
\]
and so $ \mch(A_q) = \N$.  Thus, we consider only decreasing sequences 
that are asymptotically strictly decreasing. 

For $x \in X$, the \emph{representation function} $r_{A,h}(x)$ counts 
the number of $h$-tuples 
$(a_1,\ldots, a_h) \in A \times \cdots \times A$ such that
$x = a_1 + \cdots + a_h$.  The representation function determines the sumset: 
 $hA = \{x\in X: r_{A,h}(x) > 0\}$. 
 In the additive semigroup \Z, we have $r_{\Z,h}(x) = \infty$ for all $x \in \Z$ and $h \geq 2$.

\bt[Nathanson~\cite{nath26Aus} ]            \label{intersect:theorem:finiteReps}
Let $h \in \N$.  Let $X $ be an additive abelian semigroup such that $r_{X ,h}(x) < \infty$ 
for all $x \in X $.  Let $A$ be a subset of $X $ 
and let $(A_q)_{q=1}^{\infty}$ be an asymptotically strictly decreasing sequence 
of subsets of $X$ such that $A = \bigcap_{q=1}^{\infty} A_q$. 
Then 
\[
hA = \bigcap_{q=1}^{\infty} hA_q 
\] 
and $h \in \mch(A_q)$.  
If $r_{X ,h}(x) < \infty$ for all $x \in X$ and $h \geq 2$, then 
\[
hA = \bigcap_{q=1}^{\infty} hA_q 
\]
for all $h \in \N$ and $\mch(A_q) = \N$.  
\et

\bt                    \label{intersect:theorem:nonnegative}
If $A$ is a set of nonnegative integers and if $(A_q)_{q=1}^{\infty}$ is a
 strictly decreasing sequence 
of sets of nonnegative integers such that $A = \bigcap_{q=1}^{\infty} A_q$, then 
$\mch(A_q) = \N$.
\et

\begin{proof}
The additive semigroup $\N_0$ satisfies $r_{\N_0 ,h}(x) < \infty$  
for all $h \in \N$ and $x \in \N_0$. 
\end{proof}

A sequence $(A_q)_{q=1}^{\infty}$  of sets of integers is \emph{uniformly bounded below} 
if there is an integer $m_0$ 
such that $\min A_q \geq m_0$ for all $q \in \N$.

\bt                   \label{intersect:theorem:BoundedBelow}
If $(A_q)$ is a strictly decreasing sequence 
of sets of integers that are uniformly bounded below, then $\mch(A_q) = \N$.  
\et 

\begin{proof}
Because the sequence $(A_q)_{q=1}^{\infty}$ is uniformly bounded below,
there is an integer $m_0$ such that $x \geq m_0$ for all $q \in \N$ and  $x \in A_q$.
Let $A = \bigcap_{q=1}^{\infty} A_q$.  For all $h \in \N$, 
if $x \in \bigcap_{q=1}^{\infty} hA_q$, then, for all $q \in \N$, 
there is an $h$-tuple $(a_{1,q},\ldots, a_{h,q}) \in A_q^h$ such that 
$x = a_{1,q} + \cdots + a_{h,q}$.  For all $j = 1, \ldots, h$, we have 
\[
m_0 \leq a_{j,q}= x - \sum_{\substack{i=1\\ i \neq j}}^h a_{i,q}\leq x - (h-1)m_0
\]
There are only finitely many $h$-tuples of integers $(a'_1,\ldots, a'_h)$ such that 
\[
m_0 \leq a'_j   \leq x - (h-1)m_0
\]
for all $j = 1,\ldots, h$, and so there is an $h$-tuple $(a'_1,\ldots, a'_h)$ such that 
\[
(a'_1,\ldots, a'_h) = (a_{1,q},\ldots, a_{h,q}) 
\] 
for infinitely many $q$.  For all $j = 1,\ldots, h$, we have 
$a'_j \in A_q$ for infinitely many $q$.  Because the sets $A_q$ are decreasing, 
we have $a'_j \in A_q$ for all $q \in \N$ and so $a'_j \in A = \bigcap_{q=1}^{\infty} A_q$. 
Therefore, $(a'_1,\ldots, a'_h) \in A^h$ and $x = a'_1+ \cdots + a'_h \in hA$. 
It follows that $hA =  \bigcap_{q=1}^{\infty} hA_q$ for all $h \in \N$. 
\end{proof} 

The proof of Theorem~\ref{intersect:theorem:BoundedBelow}
 is essentially the proof of 
Theorem~\ref{intersect:theorem:finiteReps} in~\cite{nath26Aus}.

\section{Arithmetical structure}

Let $a+m\Z$ be the congruence class of $a$ modulo $m$ and, for $u,v \in \Z$, let 
\[
[u,v]+m\Z = \bigcup_{a=u}^v (a+m\Z).
\]

\bt            \label{intersect:theorem:1} 
For every  integer $h_0 \geq 3$,  there is a strictly 
decreasing sequence $(A_q)_{q=1}^{\infty}$ of sets of  integers such that 
\[
\mch(A_q) = \{1\} \cup \{h_0, h_0+1,h_0+2,\ldots \}.
\]
\et

\begin{proof}
Let $h_0 \geq 3$.   
We shall construct  a strictly 
decreasing sequence $(A_q)_{q=1}^{\infty}$ of sets of  integers with 
\beq                   \label{intersect:PS1-1}
A = \bigcap_{q=1}^{\infty} A_q
\eeq
 such that 
\beq                   \label{intersect:PS1-3}
hA = \bigcap_{q=1}^{\infty} hA_q   
\eeq 
for all $h  \geq h_0$, but 
\beq                   \label{intersect:PS1-2}
h A \neq  \bigcap_{q=1}^{\infty}  h A_q 
\eeq
for all $h \in \{2,  \ldots, h_0-1 \}$.

Let $s$ be a positive integer and let 
\[
m = (h_0-1)s+2 
\]
and 
\[
A = [0,s]+m\Z. 
\]
For all $h \in \N$, we have 
\[
h A = [0,h s] + m\Z. 
\] 
Also, because $0 \in A$, the sequence 
of sumsets is increasing:   
 $h A \subseteq (h+1)A$ for all $h \in \N$.  

For $h \in \{1,2,  \ldots, h_0-1 \}$, we have 
\[
h A \subseteq (h_0-1)A = [0,(h_0-1)s] + m\Z = [0,m-2] + m\Z. 
\] 
Therefore,  
\beq            \label{intersect:1} 
(m-1+m\Z) \cap h A = \emptyset \qquad \text{for all $h \in \{1,2,  \ldots, h_0-1 \}$.}
\eeq
Because $h_0s \geq  (h_0-1)s+1 = m-1$, for all $h \geq h_0$ we have 
\[
\Z =  [0,m-1] + m\Z = [0,h_0s] + m\Z  =  h_0A \subseteq h A \subseteq \Z
\] 
and so 
\beq            \label{intersect:3} 
h A =\Z    \qquad\text{for all $h \geq h_0$.} 
\eeq

For all $q \in \N$, let 
\[
B_q =  \{ m-1+ mr: r \geq q\} 
\]
and 
\[
A_q = A \cup B_q.
\] 
The sequences $(B_q)_{q=1}^{\infty}$ and  $(A_q)_{q=1}^{\infty}$ are strictly decreasing with  
\[
\bigcap_{q=1}^{\infty} B_q = \emptyset 
\]
and 
\[
\bigcap_{q=1}^{\infty} A_q = \bigcap_{q=1}^{\infty} (A \cup B_q) = A.
\] 
Moreover, 
\[
B_q \subseteq m-1 + m\Z. 
\]
Relation~\eqref{intersect:1} implies  
\[
B_q \cap h A = \emptyset 
\]
for all $q \in \N$ and $h \in \{1,2,  \ldots, h_0-1 \}$. 

We have 
\[
 [0,m-2] + m\Z  =  [0,(h_0-1)s] + m\Z = (h_0-1)A \subseteq (h_0-1)A_q.
\]
Because $m\Z \subseteq A$ and $m-1+mq \in B_q$, for all $q \in \N$ and 
$h \geq 2$, we have  
\begin{align*} 
m-1+m\Z & = (m-1+mq) + m\Z  \subseteq B_q  + A \\ 
& \subseteq 2(A \cup B_q ) = 2A_q \subseteq h A_q 
\end{align*} 
and so 
\beq           \label{intersect:2}
m-1+m\Z \subseteq \bigcap_{q=1}^{\infty} h A_q \qquad\text{for all $h \geq 2$.}
\eeq
Comparing~\eqref{intersect:1} and~\eqref{intersect:2}, we see that 
\[ 
h A \neq \bigcap_{q=1}^{\infty} h A_q 
\] 
and so $h \notin \mch(A_q)$ for all $h \in \{ 2,\ldots, h_0-1 \}$. 

For all $h \geq h_0-1 \geq 2$ and $q \in \N$, we have 
\begin{align*} 
\Z & =  [0,m-1]+m\Z   = \left(  [0,m-2]+m\Z \right) \cup \left(m-1+m\Z   \right) \\ 
& = \left( [0,(h_0-1)s] + m\Z \right) \cup \left(m-1+m\Z   \right) \\ 
& \subseteq  (h_0-1)A \cup 2A_q  \\
& \subseteq (h_0-1)A_q \\
& \subseteq  h A_q  \\
& \subseteq \Z  
\end{align*} 
and so 
\beq            \label{intersect:4} 
\bigcap_{q=1}^{\infty} h A_q  =\Z    \qquad\text{for all $h \geq h_0-1$.} 
\eeq
Comparing~\eqref{intersect:3} and~\eqref{intersect:4}, we see that 
\[
h A =\bigcap_{q=1}^{\infty} h A_q    \qquad\text{for all $h \geq h_0$}
\]
and so $h \in \mch(A_q)$ for all $h \geq h_0$. 
This completes the proof. 
\end{proof}

\bt                   \label{intersect:theorem:2}
For every  integer $d \geq 2$, there exists a strictly decreasing sequence 
$(A_q)_{q=1}^{\infty}$ of sets of integers with 
\[
\bigcap_{q=1}^{\infty} A_q = \emptyset 
\]
such that  
\[
\bigcap_{q=1}^{\infty} hA_q= \emptyset 
\]
if and only if $d$ does not divide $h$. 
\et

\bc 
For every  integer $d \geq 2$, there exists a strictly decreasing sequence $(A_q)_{q=1}^{\infty}$ 
of sets of integers such that 
\[
\mch(A_q) = \{h \in \N: h \not\equiv 0 \pmod{d}\}.
\]
 \ec

\begin{proof}
Let $(g_j)_{j=1}^{\infty}$ be a strictly increasing sequence of positive integers and let 
\[
G_q= \prod_{j=1}^q g_j
\]
for all $q \in \N$.   The sequence $(G_q)_{q=1}^{\infty}$ is strictly increasing 
and defines an additive system (Nathanson~\cite{nath14a}).  
The integer $G_q/G_k$ is divisible by $g_{k+1}$ for all $q \geq k+1$. 

For all $q \in \N$, let 
\[
A_q = \{G_r, -(d-1)G_r: r \geq q\}.
\]
The sequence of sets $(A_q)_{q=1}^{\infty}$ is strictly decreasing and 
\[
A = \bigcap_{q=1}^{\infty} A_q = \emptyset.  
\]
If $h = d'd$, then $h = d'(d-1)+d'$.  For all $q \in \N$, we have $G_q \in A_q$ and 
\begin{align*}
0 & = d'  \left( (d-1) G_q -(d-1)G_q \right) \\ 
& = d'(d-1) G_q + d'(-(d-1)G_q) \\ 
& \in hA_q.
\end{align*}
It follows that if $d$ divides $h$, then $0 \in \bigcap_{q=1}^{\infty} hA_q$. 

Conversely, if $0 \in \bigcap_{q=1}^{\infty} hA_q$, then, for all $q \in \N$, 
we have $0 \in hA_q$ and so there exist sequences of nonnegative integers  
$(x_r)_{r=q}^{\infty}$ and $(y_r)_{r=q}^{\infty}$ such that 
\[
h = \sum_{r=q}^{\infty} (x_r+y_r)  
\]
and 
\begin{align*} 
0 & = \sum_{r=q}^{\infty} (x_rG_r + y_r(-(d-1)G_r)) \\ 
& = \sum_{r=q}^{\infty} (x_r - (d-1) y_r) G_r \\ 
&  = \sum_{r=q}^{\infty} z_rG_r
\end{align*} 
where $z_r = x_r - (d-1) y_r$ for all $r \geq q$, and 
\[
|z_r| = |x_r - (d-1) y_r| \leq (d-1)(x_r+y_r) \leq (d-1)h.
\]

Suppose that $z_r \neq 0$ for some $r \geq q$.   Let $k$ be the smallest integer such that 
$k \geq q$ and $z_{k}\neq 0$.  Then 
\[
z_{k} G_{k} = - \sum_{r=k+1}^{\infty} z_rG_r. 
\] 
Because  $g_{k+1}$ divides $G_r/G_k$ for all $r \geq k+1$, 
there is a nonzero integer $\ell$ such that 
\[
0 \neq z_{k} = - \sum_{r=k+1}^{\infty} z_r  \left( \frac{G_r}{G_{k} } \right) = g_{k+1} \ell.
\]
This implies 
\[
g_q  \leq g_k < g_{k+1} \leq |z_k| \leq (d-1)h
\]
for all $q$, which is absurd because $\lim_{j\rightarrow \infty} g_j = \infty$. 
It follows that, for all $r \geq q$, we have $z_r = x_r - (d-1) y_r = 0$ and so 
\[
h = \sum_{r=q}^{\infty} (x_r+y_r) = \sum_{r=q}^{\infty} ((d-1)y_r +y_r) 
 = d\sum_{r=q}^{\infty}  y_r \equiv 0 \pmod{d}. 
\]
This proves that $0 \in \bigcap_{q=1}^{\infty} hA_q$ if and only if 
$h \equiv 0 \pmod{d}$. 

If $n \in  \bigcap_{q=1}^{\infty} hA_q$ and $n \neq 0$, then, for all $q \in \N$, 
we have $n \in hA_q$ and there is a sequence $(z_r)_{r=q}^{\infty} $  
of  integers such that 
\[
n = \sum_{r=q}^{\infty} z_rG_r = G_q \sum_{r=q}^{\infty} z_r \left(\frac{G_r}{G_q} \right) 
\]
and so $|n| \geq G_q$ for all $q$, which is absurd 
because $\lim_{j\rightarrow \infty} G_q = \infty$.  Therefore, 
\[
\bigcap_{q=1}^{\infty} hA_q = 
\begin{cases}
\{0\} & \text{if $h \equiv 0 \pmod{d}$}\\
\emptyset & \text{if $h \not\equiv 0 \pmod{d}$.}
\end{cases}
\]
Because $hA = \emptyset$ for all $h \in \N$, we have  
$\mch(A_q) = \{h \in \N: h \not\equiv 0 \pmod{d}\}$.  
This completes the proof. 
\end{proof}

From Theorems~\ref{intersect:theorem:1} and~\ref{intersect:theorem:2} we obtain 

\bc
For every $h_0 \geq 2$ 
there exist sequences  $(A_q)_{q=1}^{\infty}$ such that 
$\{1,\ldots, h_0-1\} \subseteq \mch(A_q)$ but $h_0 \notin \mch(A_q)$ 
and also that there exist sequences  $(A_q)_{q=1}^{\infty}$ such that 
$\{1,  h_0 \} \subseteq \mch(A_q)$ but $\{2,3, \ldots, h_0-1\} \cap \mch(A_q) = \emptyset$.
\ec

\end{document}